\RequirePackage{ifpdf}
\ifpdf 
\documentclass[pdftex]{sigma}
\else
\documentclass{sigma}
\fi

\usepackage{mathrsfs,dsfont}
\usepackage{stmaryrd}

\newcommand{\R}{\mathbb{R}}
\newcommand{\pp}[2]{\frac{\partial #1}{\partial #2}}
\newcommand{\com}[2]{[\,#1,#2\,]}

\newcommand{\n}{^{(n)}}
\newcommand{\ii}{^{(\infty )}}

\newcommand{\vv}{\mathbf{v}}
\newcommand{\ww}{\mathbf{w}}

\newcommand{\bt}{\pmb{\tau}}

\newcommand{\bti}{\bt\ii}

\newcommand{\mui}{\mu\ii}

\newcommand{\g}{\mathfrak{g}}
\newcommand{\Dt}{\mathbb{D}}

\newcommand{\D}{\mathcal{D}}
\newcommand{\Dn}{\D\n}
\newcommand{\Di}{\D\ii}
\newcommand{\G}{\mathcal{G}}
\newcommand{\Gn}{\G\n}
\newcommand{\Gi}{\G\ii}
\newcommand{\J}{J}
\newcommand{\Jn}{\J^n}
\newcommand{\Ji}{\J^\infty}
\newcommand{\Z}{Z}
\newcommand{\Zn}{\Z\n}

\newcommand{\sminus}{\mkern-2mu\setminus\mkern-2mu}

\def\comp{\raise 1pt \hbox{$\,\scriptstyle\circ\,$}}

\def\dom{\mathop{\rm dom}\nolimits}
\def\dim{\mathop{\rm dim}\nolimits}

\def\subs #1#2{#1_1,\ldots,#1_{#2}}

\def\sups #1#2{#1^1,\ldots,#1^{#2}}
\def\psups #1#2{(#1^1,\ldots,#1^{#2})}

\def\Pa#1{\left (\,{#1}\,\right )}
\def\bpa#1{\bigl(\,{#1}\,\bigr)}
\def\br#1{[\mskip2mu{#1}\mskip 2mu]}

 \def\rbox#1{\hbox{\rm #1}}

\def\roq#1{\qquad \rbox{#1}\qquad }

\def\mcases#1{\left\{\enspace\seq{#1}\right.}
\def\seq#1{{\def\\{\cr}\vcenter{\openup 1\jot\halign{$\displaystyle
   ##\hfil$&&\hskip20pt$\displaystyle##\hfil$\cr #1\cr}}}}

\numberwithin{equation}{section}

\newtheorem{Theorem}{Theorem}[section]
{\theoremstyle{definition}
\newtheorem{Definition}[Theorem]{Definition}
\newtheorem{Remark}[Theorem]{Remark}
\newtheorem{Example}[Theorem]{Example}}

\begin{document}

\allowdisplaybreaks

\renewcommand{\thefootnote}{$\star$}

\renewcommand{\PaperNumber}{077}

\FirstPageHeading

\ShortArticleName{On the Structure of Lie Pseudo-Groups}

\ArticleName{On the Structure of Lie Pseudo-Groups\footnote{This paper is a
contribution to the Special Issue ``\'Elie Cartan and Dif\/ferential Geometry''. The
full collection is available at
\href{http://www.emis.de/journals/SIGMA/Cartan.html}{http://www.emis.de/journals/SIGMA/Cartan.html}}}

\Author{Peter J. OLVER~$^\dag$, Juha POHJANPELTO~$^\ddag$ and Francis VALIQUETTE~$^\dag$}

\AuthorNameForHeading{P. Olver, J. Pohjanpelto and F. Valiquette}

\Address{$^\dag$~School of Mathematics, University of Minnesota, Minneapolis, MN 55455, USA} 
\EmailD{\href{mailto:olver@math.umn.edu}{olver@math.umn.edu}, \href{mailto:valiq001@math.umn.edu}{valiq001@math.umn.edu}}
\URLaddressD{\url{http://www.math.umn.edu/~olver/}, \url{http://www.math.umn.edu/~valiq001/}}

\Address{$^\ddag$~Department of Mathematics, Oregon State University, Corvallis, OR 97331, USA}
\EmailD{\href{mailto:juha@math.oregonstate.edu}{juha@math.oregonstate.edu}}
\URLaddressD{\url{http://oregonstate.edu/~pohjanpp/}}

\ArticleDates{Received March 31, 2009, in f\/inal form July 08, 2009;  Published online July 23, 2009}

\Abstract{We compare and contrast two approaches to the structure theory for Lie pseudo-groups, the f\/irst due to Cartan, and the second due to the f\/irst two authors.  We argue that the latter approach of\/fers certain advantages from both a theoretical and practical standpoint.}


\Classification{58A15; 58H05}

\renewcommand{\thefootnote}{\arabic{footnote}}
\setcounter{footnote}{0}

\section{Introduction}

The aim of this paper is to compare the structure theory for Lie pseudo-groups developed by the f\/irst two authors in \cite{OP-2005}  with the classical Cartan theory, \cite{C-1904-0,C-1904,C-1937}.  The former relies on the contact structure of the inf\/inite dif\/feomorphism jet bundle, whereas Cartan's is based on the prolongation of exterior dif\/ferential systems.  We show that the two theories are isomorphic in the case of transitive Lie pseudo-groups, but lead to dif\/ferent structure equations when dealing with intransitive pseudo-group actions.  We then argue that the former theory of\/fers some distinct advantages over the Cartan structure theory in this situation.

Our reference point is the well-established structure theory for f\/inite-dimensional Lie groups.  Let $G$ be a Lie group of dimension $r = \dim G$. The commutators
\begin{gather}\label{g}
\com{\vv_j}{\vv_k} = \sum_{i=1}^r   C^i_{jk}\vv_i
\end{gather}
between the inf\/initesimal generators -- that is, a basis $\subs \vv r$ for its Lie algebra $\g$ -- prescribe the \emph{structure constants} $C^i_{jk} = -  C^i_{kj}$, which serve to uniquely characterize a connected Lie group~$G$ up to a discrete subgroup.  Equivalently, the structure of $G$ can be based on the \emph{Maurer--Cartan structure equations}
\begin{gather}\label{MC}
d \mu^i = -\sum_{j<k}   C^i_{jk} \mu^j \wedge \mu^k,
\end{gather}
satisf\/ied by the Maurer--Cartan one-forms $\sups \mu r$, which form the dual basis of the dual space $\g^*$.  It is noteworthy that the \emph{same} structure constants appear in both the Maurer--Cartan structure equations and the commutator relations.

The key obstruction hindering an immediate generalization of the f\/inite-dimen\-sional structure theory to inf\/inite-dimensional Lie pseudo-groups is the lack of an appropriate abstract object to represent the pseudo-group itself.  Thus, at least in our current state of knowledge, Lie pseudo-groups are inextricably bound to the manifold on which they act.  The appropriate Maurer--Cartan forms thus must be suitably invariant dif\/ferential forms living on the manifold or, better, on some bundle associated with the action.  The approach developed in~\cite{OP-2005} is based on the bundle of inf\/inite order jets of pseudo-group transformations and the invariant contact forms thereon.

More specif\/ically, the bundle of inf\/inite jets of local dif\/feomorphisms belonging to the pseudo-group forms a subbundle -- indeed a subgroupoid -- of the inf\/inite dif\/feomorphism jet bundle.  The Maurer--Cartan forms will be identif\/ied as the right-invariant\footnote{Alternatively, one can use the left-invariant forms.  As in the references, for specif\/icity, we focus on the right-invariant constructions here.} contact forms on the latter jet bundle, and their structure equations are readily found.   Restricting the dif\/feomorphism-invariant contact forms to the pseudo-group jet subbundle results in a system of Maurer--Cartan forms for the pseudo-group, whose structure equations are obtained by restriction of the dif\/feomorphism structure equations.  Remarkably, the restricted invariant contact forms, which are no longer linearly independent, satisfy a collection of linear algebraic constraints that can be immediately obtained by lifting the inf\/initesimal determining equations for the pseudo-group.  This allows us to immediately establish a complete system of structure equations for \emph{any} Lie pseudo-group directly from its inf\/initesimal determining equations, thereby avoiding the more cumbersome and unintuitive prolongation construction advocated by Cartan.   We emphasize that the method \emph{does not} rely on the explicit formulas for the Maurer--Cartan forms, and only needs elementary linear algebra to extract the complete structure equations.  Moreover, the construction can be readily implemented in any coordinate system on the underlying manifold, and avoids having to identify the invariants and work in the specially adapted coordinates as required by Cartan's method. As a result, the theory can be immediately applied in a broad range of examples, and the necessary algorithms are straightforwardly implemented using standard symbolic software packages such as {\sc Mathematica} or {\sc Maple}.

Another advantage of the contact form approach is that it applies equally well to both transitive and intransitive pseudo-groups, and naturally includes f\/inite-dimensional Lie transformation groups as a particular case.  In the transitive case, we show that the Cartan structure equations are isomorphic to those satisf\/ied by the invariant contact forms.  However, a direct isomorphism is no longer valid in the more challenging case of intransitive pseudo-group actions.  Furthermore, the Maurer--Cartan structure equations established here are directly dual to the commutator equations for the inf\/initesimal generators of the pseudo-group, and, moreover, coincide with the structure equations \eqref{MC} when the pseudo-group is of f\/inite type.  This is in contrast to Cartan's version, which, in particular, produces nonzero structure constants/functions for intransitive actions of abelian pseudo-group and Lie group actions, \cite{F-2008}, thus making the connections between the structure equations and the Lie algebra of inf\/initesimal generators somewhat obscure.

One of the main results of Cartan is that any Lie pseudo-group, after a f\/inite number of prolongations, is characterized by leaving a coframe and a certain number of functions invariant.  By virtue of the Cartan--K\"ahler theorem, \cite{BCGGG-1991,O-1995}, Cartan's structure equations serve as integrability conditions on the invariant coframe.   The invariant coframe constructed by Cartan depends on the realization of the pseudo-group action, and two pseudo-groups that are isomorphic  in the sense of Cartan, \cite{C-1904,C-1937,S-2000}, can have non-isomorphic Cartan structure equations.  On the other hand, our Maurer--Cartan structure equations are always isomorphic under Cartan's notion of isomorphism of pseudo-groups, \cite{F-2008,F-2009}.

Pertinent references on the general theory of Lie pseudo-groups include the classical works of Lie, \cite{L-1883,L-1891}, Cartan, \cite{C-1904-0,C-1904,C-1937}, and Vessiot, \cite{V-1903}, along with a variety of contemporary treatments, \cite{Chrastina,Ehresmann,K-1989,K-1959,K-1961,LisleReidi,OP-2005,P-1978,SS-1965,S-2000}.  The basics of jet bundles, contact forms, and the variational bicomplex can be found, for instance, in  \cite{A-1989,O-1995}.  Applications of these results in the method of moving frames for pseudo-groups can be found in \cite{OP-2008,OP-2009}.

\section{The dif\/feomorphism pseudo-group}\label{section diff}

We begin by describing the structure of the most basic pseudo-group.
Let $M$ be a smooth $m$-dimensional manifold and write $\D=\D(M)$ for the pseudo-group of all local\footnote{Our notational conventions allow the domain of def\/inition of a map $\varphi \colon  M \to M$ to be a proper open subset: $\dom \varphi \subset M$.  Also, when we write $Z = \varphi(z)$ we implicitly assume $z \in \dom \varphi$.} dif\/feomorphisms $\varphi\colon  M \to M$.  For each $0 \leq n \leq \infty $, let $\Dn \subset \Jn(M,M)$ denote the bundle of their $n$-th order jets.  We remark that $\Dn$ carries the structure of a groupoid, \cite{Mackenzie}, whose multiplication is provided by algebraic composition of Taylor series (when def\/ined).  There are natural right and left actions of $\D$ on the jet bundles $\Dn$, denoted by $R_ \varphi $ and $L_ \varphi $, respectively.

Local coordinates $(z,\Zn)$ on $\Dn$ are provided by
a system of source coordinates $z=(z^1,\ldots,z^m)$ on $M$, target coordinates $Z = (Z^1,\ldots,Z^m)$ also on $M$, and jet coordinates $Z^b_A$ representing the partial derivatives $\partial ^{k}\varphi^b(z)/\partial z^{a_1} \cdots \partial z^{a_k}$ of the local dif\/feomorphism $Z =\varphi (z)$.  Here $A=(a_1,\dots, a_k)$, with $1 \leq a_\nu \leq m$, indicates a multi-index of order $k = \#A \leq n$.
In what follows, we will consistently use lower case letters, $z,x,u,\ldots $
for the source coordinates and the corresponding upper case letters $Z,X,U,\ldots $ for the target coordinates.

At inf\/inite order, the cotangent bundle $T^\star \Di \subset T^\star \Ji(M,M)$ naturally splits into horizontal and vertical (contact) components, spanned respectively by the coordinate dif\/ferentials $dz^1,\ldots,dz^m$, and the basic {\it contact forms}
\begin{gather}\label{cf}
\Upsilon_A^b = dZ_A^b - \sum_{a=1}^m  Z_{A,a}^b dz^a,\qquad b=1,\ldots,m,\qquad \#A \geq 0.
\end{gather}
The decomposition of $T^\star \Di$  accordingly splits the dif\/ferential
$d=d_M + d_G$, where the subscript on the vertical dif\/ferential $d_G$ refers to the groupoid structure of $\Di$.  In particular, if $F(z,\Zn)$ is any dif\/ferential function, then
\[
d_M F=\sum_{a=1}^m  (\Dt_{z^a} F) dz^a,\qquad d_G F = \sum_{b=1}^m \sum_{\# A \geq 0} \pp{F}{Z^b_A} \Upsilon^b_A,
\]
where
\begin{gather}\label{Dt}
\Dt_{z^a} = \pp{}{z^a}+\sum_{b=1}^m \sum_{\# A \geq 0} Z^b_{A,a}\pp{}{Z^b_A},\qquad a=1,\ldots,m,
\end{gather}
denotes the coordinate \emph{total derivative operators}.

Since the target coordinate functions $Z^a\colon \Di \to \R$ are clearly invariant under the right action of $\D$, so are their dif\/ferentials $dZ^a$.  The splitting of the dif\/ferential into horizontal and contact components is also right-invariant.
This implies that the one-forms
\begin{gather}\label{invariant horizontal coframe}
\sigma^a = d_M Z^a = \sum_{b=1}^m  Z^a_b dz^b,\qquad a=1,\ldots,m,
\end{gather}
form an invariant horizontal coframe, while
\begin{gather}\label{zero mc forms}
\mu^a=d_G Z^a = \Upsilon^a = dZ^a - \sum_{b=1}^m  Z^a_b dz^b,\qquad a=1,\ldots,m,
\end{gather}
are the zero-th order invariant contact forms.  Writing the horizontal component of the dif\/ferential of a dif\/ferential function $ F \colon \Di \to \R$ in terms of the invariant horizontal coframe~\eqref{invariant horizontal coframe},
\[
d_M F = \sum_{a=1}^m  (\Dt_{Z^a} F) \sigma^a,
\]
serves to def\/ine the dual invariant total dif\/ferential operators
\begin{gather}\label{differential operators}
\Dt_{Z^a}=\sum_{b=1}^m  w_a^b \,\Dt_{z^b},\qquad a=1,\ldots,m,
\end{gather}
where
\[
\bigl( w_a^b\big(z,Z^{(1)}\big) \bigr)=\left(\pp{Z^b}{z^a}\right)^{-1}
\]
denotes the inverse Jacobian matrix.
Thus, higher-order right-invariant contact forms are obtained by successively applying the invariant dif\/ferential operators  \eqref{differential operators} to the zero-th order invariant contact forms \eqref{zero mc forms}:
\begin{gather}\label{mc forms}
\mu^a_A=\Dt_Z^A \mu^a= \Dt_{Z^{a_1}}\cdots\Dt_{Z^{a_k}}\mu^a,\qquad a=1,\ldots,m,\qquad \# A \geq 0.
\end{gather}
The dif\/ferential operators $\Dt_{Z^1},\ldots\!,\Dt_{Z^m}$ mutually commute, so the order of dif\/ferentiation in~\eqref{mc forms} is immaterial.  As in \cite{OP-2005}, we interpret the right-invariant contact forms $\mui = ( \ldots \mu^a_A  \ldots )$  as the \emph{Maurer--Cartan forms} for the dif\/feomorphism pseudo-group $\D$, and they, together with the horizontal forms
\eqref{invariant horizontal coframe} provide a right-invariant coframe on $\D\ii$.

The dif\/feomorphism structure equations satisf\/ied by the Maurer--Cartan forms are easily established, \cite{OP-2005}.  They can be concisely expressed by introducing the vector-valued Maurer--Cartan formal power series $\mu\llbracket H \rrbracket = (\mu^1\llbracket H \rrbracket , \ldots, \mu^m\llbracket H \rrbracket )^T$,  with components
\begin{gather}\label{mc power series}
\mu^a \llbracket H \rrbracket = \sum_{\# A \geq 0} \frac{1}{A!}\,\mu_A^a \,H^A,\qquad a=1,\ldots,m.
\end{gather}
Here $H = (H^1,\ldots,H^m)$ are formal power series parameters, while
$A! = i_1!\,i_2!\cdots i_m!$,  where $i_l$ stands for the number of occurrences of the integer $l$ in $A$.
The structure equations for the right-invariant forms $\mu^a_A$ are obtained by comparing the coef\/f\/icients of the various powers of $H$ in the power series identity
\begin{gather}\label{diff structure eq1}
d\mu \llbracket H \rrbracket = \nabla \mu \llbracket H \rrbracket \wedge (\mu \llbracket H \rrbracket - dZ ),
\end{gather}
where $dZ = (dZ^1,\ldots,dZ^m)^T$, and where
$
\nabla \mu \llbracket H \rrbracket = \left( {\partial \mu^a \llbracket H \rrbracket}/{\partial H^b}\right)
$
denotes the $m \times m$ Jacobian matrix obtained by formal dif\/ferentiation of the power series \eqref{mc power series} with respect to the parameters.
The complete structure equations for the dif\/feomorphism pseudo-group are then furnished by equations \eqref{diff structure eq1} together with the equations
 \begin{gather}\label{diff structure eq2}
d\sigma = \nabla \mu \llbracket 0 \rrbracket \wedge \sigma
\end{gather}
for the invariant horizontal forms $\sigma=(\sigma^1,\ldots,\sigma^m)^T$.
We restrict the structure equations \eqref{diff structure eq1} to a target f\/iber $(\bti)^{-1}(Z)\subset \Di$  to obtain the \emph{Maurer--Cartan structure equations} for the dif\/feomorphism pseudo-group.
This amounts to setting
\begin{gather}\label{key equality}
0=dZ=\sigma+\mu\llbracket 0 \rrbracket,\roq{so that} \sigma=-\mu\llbracket 0 \rrbracket.
\end{gather}
Consequently, the structure equations \eqref{diff structure eq2} for the horizontal forms $\sigma$ become identical with the structure equations for the zero-th order Maurer--Cartan forms $\mu^a =\mu^a\llbracket 0 \rrbracket$.

\begin{Theorem}\label{MCseqs}
The Maurer--Cartan structure equations for the diffeomorphism pseudo-group $\D$ are
\begin{gather}\label{diff structure eq fiber}
d\mu \llbracket H \rrbracket = \nabla \mu \llbracket H \rrbracket \wedge \mu \llbracket H \rrbracket.
\end{gather}
\end{Theorem}

\begin{Example}
For the pseudo-group $\D(\R)$ of local dif\/feomorphisms of $M = \mathbb{R}$, the Maurer--Cartan power series is
\[
\mu\llbracket H \rrbracket = \sum_{n=0}^\infty  \mu_n\, \frac{H^n}{n!},
\]
where $\mu_n=\Dt_X^n\mu_0$, $n = 0, 1,2 \ldots$, are the right-invariant contact forms on $\Di$.  The individual components of \eqref{diff structure eq fiber} yield the expressions
\begin{gather*}
d\mu_n= \sum_{i=0}^n \binom{n}{i}\mu_{i+1} \wedge \mu_{n-i}
= -\sum_{i=0}^{[(n+1)/2]} \frac{n-2\,i+1}{n+1}\binom{n+1}{i}\mu_i \wedge\mu_{n+1-i},\qquad n\geq 0,
\end{gather*}
which reproduce the structure equations found by Cartan, \cite[equation~(48)]{C-1904}.
\end{Example}

Expanding the power series \eqref{mc power series}, we f\/ind that  the Maurer--Cartan structure equations \eqref{diff structure eq fiber} have the individual components
\begin{gather}\label{diffeomorphism MCx}
d\mu^a_C=\sum_{C = (A,B)}\sum_{b=1}^m  {C\choose A}\,\mu^a_{A,b}\wedge\mu^b_B,
\end{gather}
involving the multinomial coef\/f\/icients
\[
 {C\choose A} = \frac{C\:!}{A\:!\,B\:!} \qquad \text{when} \qquad C = (A,B) = (a_1,\ldots ,a_k,b_1,\ldots ,b_l)
 \]
 is the union of two multi-indices $A= (a_1,\ldots ,a_k)$ and $B=(b_1,\ldots ,b_l)$, either of  which can be empty.

\begin{Remark}
Since the higher order Maurer--Cartan forms $\mu^a_{C}$ are def\/ined by \eqref{mc forms}, their structure equations \eqref{diffeomorphism MCx} can also be derived by Lie dif\/ferentiating the structure equations for the zero-th order invariant contact forms $\mu^a$.  By direct computation
\[
d\mu^a=\sum_{b=1}^m  \mu^a_{b}\wedge \big(\mu^b-dZ^b\big),
\]
and, from the Leibniz rule, we obtain
\begin{gather*}
d\mu^a_{C}=d\left(\mathbb{D}^C_Z \mu^a\right)=\mathbb{D}^C_Z \left(d\mu^a\right)= \mathbb{D}^C_Z\left( \sum_{b=1}^m \mu^a_{b}\wedge (\mu^b-dZ^b) \right) \\
\phantom{d\mu^a_{C}=d\left(\mathbb{D}^C_Z \mu^a\right)=\mathbb{D}^C_Z \left(d\mu^a\right)}{}
= \sum_{C=(A,B)}\sum_{b=1}^m \binom{C}{A}\,\mu^a_{A,b}\wedge\left(\mu^b_{B}-d\big(\mathbb{D}_Z^B Z^b\big)\right).
\end{gather*}
The last term, $d(\mathbb{D}_Z^B Z^b)$, is non-trivial only when $\#B = 0$.
Restricting the last equation to a~target f\/iber $(\bti)^{-1}(Z)$ we recover the Maurer--Cartan structure equations \eqref{diffeomorphism MCx}.
\end{Remark}

\section{Lie pseudo-groups}\label{structure pseudo-group section}

Let $\G \subset \D$ be a sub-pseudo-group acting on $M$ and let $\Gn \subset \Dn$, $0 \leq n \leq \infty $, denote the corresponding subgroupoid of $n$-th order jets of its local dif\/feomorphisms.  Roughly speaking, $\G$~is called a Lie pseudo-group provided that it can be identif\/ied as the solution space to a system of partial dif\/ferential equations. There are several variants of the precise technical requirements to be found in the literature; ours are the following.

\begin{Definition}\label{definition Lie pseudo-group}
A sub-pseudo-group $\G\subset \D$ is called a \emph{Lie pseudo-group} if there exists $n^\star \geq 1$ such that, for all f\/inite $n \geq n^\star$:
\begin{enumerate}\itemsep=0pt
\item
$\Gn \subset \D^{(n)}$ forms a smooth, embedded subbundle;
\item
$\pi_n^{n+1}:\G^{(n+1)}\to \Gn$ is a f\/ibration;
\item
if $j_n\varphi \subset \Gn$, then $\varphi\in \G$;
\item
$\Gn=\text{pr}^{(n-n^\star)}\G^{(n^\star)}$ is obtained by prolongation.
\end{enumerate}
The minimal value of $n^\star$ is called the \emph{order} of the Lie pseudo-group.
\end{Definition}

Thus, by condition 1, the pseudo-group jet subbundle $\Gn\subset \Dn$ is prescribed in local coordinates by a system of $n$-th order (typically nonlinear) partial dif\/ferential equations
\begin{gather}\label{determining equations}
F\n(z,Z\n)=0,
\end{gather}
known as the $n$-th order \emph{determining equations} for the Lie pseudo-group $\G$.  By construction, for any $n\geq n^\star$, the system \eqref{determining equations} is locally solvable, and its local solutions $Z=\varphi(z)$, by condition~3, are precisely the pseudo-group transformations.  Moreover, by condition~4, the determining equations in order $n > n^\star$ can be obtained by repeatedly applying the total
derivative operators~\eqref{Dt} to those of order $n^\star$.

Let $\g$ denote the local Lie algebra\footnote{By ``local Lie algebra'', we mean that the vector f\/ields $\vv \in \g$ may only be locally def\/ined, and that, for $\vv,\ww \in \g$, the Lie bracket $\com \vv\ww \in \g$ is only def\/ined on their common domain of def\/inition.} of inf\/initesimal generators of our pseudo-group $\G$, i.e., the set of locally def\/ined vector f\/ields whose f\/low maps belong to the pseudo-group. Let $\Jn \g \subset \Jn TM$, $0\leq n\leq\infty$, denote their $n$-jets.  Fiber coordinates on the vector f\/ield jet bundle $\Jn TM$  are given by $\zeta ^b_A$, for $1 \leq b \leq m$, $0 \leq \#A \leq n$, representing the partial derivatives $\partial ^{k}\zeta ^b(z)/\partial z^{a_1} \cdots \partial z^{a_k}$ of the components of a vector f\/ield
\begin{gather}\label{infinitesimal generator}
\vv = \sum_{a = 1}^m  \zeta ^a(z) \frac{\partial}{\partial z^a}
\end{gather}
written in local coordinate $z=(z^1,\ldots,z^m)$ on $M$.
By linearizing the $n$-th order pseudo-group determining equations \eqref{determining equations} at the $n$-jet of the identity transformation, we see that the
subbundle $\Jn\g \subset \Jn TM$ can locally be viewed as a system of linear partial dif\/ferential equations
\begin{gather}\label{linearized det system}
L\n\big(z,\zeta\n\big)=0,
\end{gather}
for the coef\/f\/icients $\zeta(z) = (\zeta^1(z),\dots,$ $\zeta^m(z))$ of the inf\/initesimal generators, called the \emph{infini\-te\-simal determining equations} of the pseudo-group.  In particular, if $\G$ arises as the symmetry pseudo-group of a system of partial dif\/ferential equations, then \eqref{linearized det system} is the involutive completion of the usual system of determining equations resulting from applying Lie's inf\/initesimal symmetry algorithm, \cite{O-1995}.

The Maurer--Cartan forms associated with the Lie pseudo-group $\G$ are obtained by pulling back the dif\/feomorphism Maurer--Cartan forms \eqref{mc forms} to the pseudo-group jet bundle \mbox{$\Gi{\subset} \Di$}.   The restricted Maurer--Cartan forms are, of course, no longer linearly independent.  However, the induced dependencies can, perhaps surprisingly, be explicitly prescribed with the aid of the inf\/initesimal determining equations, \cite{OP-2005}.

\begin{Theorem}\label{mc lin relations prop}
The complete set of linear dependencies among the right-invariant Maurer--Cartan forms $\mui$ is provided by the linear system
\begin{gather}\label{mc lin relations}
L\ii\big(Z,\mui\big)=0,
\end{gather}
obtained from the linear determining equations \eqref{linearized det system} by replacing the source variables $z^a$ by the corresponding target variables $Z^a$, and the infinitesimal generator jet coordinates $\zeta^b_A$ by the corresponding Maurer--Cartan forms $\mu^b_A$.
\end{Theorem}

The equations \eqref{mc lin relations} are called the \emph{lifted infinitesimal determining equations} for the Lie pseudo-group $\G$.  (See \cite{OP-2005} for additional details on the lifting process.)  Thus, the structure equations for our pseudo-group can simply be obtained by restricting the dif\/feomorphism structure equations~\eqref{diff structure eq fiber} to the solution space to the lifted inf\/initesimal determining equations~\eqref{mc lin relations}.

\begin{Theorem}
The Maurer--Cartan structure equations of a Lie pseudo-group $\G$ are obtained by imposing the linear relations prescribed by the lifted infinitesimal determining equations \eqref{mc lin relations} on the diffeomorphism Maurer--Cartan equations \eqref{diff structure eq fiber}:
\begin{gather}\label{pseudo-group structure eq}
\bpa{d\mu\llbracket H \rrbracket = \nabla\mu\llbracket H \rrbracket \wedge \mu \llbracket H \rrbracket}\big|_{L\ii(Z,\mui)=0}.
\end{gather}
\end{Theorem}

\begin{Remark} The motivation behind the need to restrict the Maurer--Cartan forms to a target f\/iber can be readily understood in the context of f\/inite-dimensional Lie group actions.  In this situation, $\bti\colon\Gi\to M$ will typically be a principal $G$ bundle, and, consequently, the independent Maurer--Cartan forms on $\Gi$ and their structure equations, when restricted to a~target f\/iber $(\bti)^{-1}(Z)\simeq G$ coincide with the usual Maurer--Cartan forms and their structure equations \eqref{MC}.

However, it is worth pointing out that, due to the appearance of the coordinates $z$ in the inf\/initesimal determining equations \eqref{linearized det system}, the basis of $\g^\star$ prescribed by the independent restricted invariant contact forms $\mu_A^b$ may vary from f\/iber to f\/iber as the target point $Z$ ranges over $M$.  Consequently, the structure coef\/f\/icients in the pseudo-group structure equations \eqref{pseudo-group structure eq} may very well be $Z$-dependent.  It is noteworthy that, when $\G$ is of f\/inite type and so represents the action of a f\/inite-dimensional Lie group $G$ on $M$, the resulting variable structure coef\/f\/icients $C^i_{jk}(Z)$ represent the {\it same} Lie algebra $\g$ and so are all similar, modulo a $Z$-dependent change of basis, to the usual {\it constant structure coefficients} associated with a f\/ixed basis of $\g^\star$.  In contrast, this is not necessarily the case for inf\/inite-dimensional intransitive pseudo-group actions. The non-constant invariants (under change of Maurer--Cartan basis) of the structure coef\/f\/icients are {\it essential invariants}, f\/irst exposed by Cartan, \cite{C-1904-0,S-2000}; see \cite{F-2008,F-2009} for further details, comparing Cartan's approach with ours. The existence of essential invariants is one of the key obstacles preventing the construction of a suitable abstract object representing the pseudo-group.
\end{Remark}

\section{Cartan structure equations}\label{cartan structure equations section}

In this section we provide a brief overview of Cartan's method for constructing the structure equations for a Lie pseudo-group.  For more detailed accounts, we refer the reader to Cartan's original works \cite{C-1904,C-1937}, and to the expository texts \cite{K-1989,S-2000}.

Given a Lie pseudo-group $\G$ acting on $M$, we choose local coordinates
\[
z = (x,y) = (\sups xs,\sups yt) = \psups zm,\qquad \quad s+t=m=\dim M,
\]
so that the pseudo-group action takes the form
\begin{gather}\label{cartan pg}
X^i=x^i,\qquad Y^\alpha=f^\alpha(x,y),\qquad i=1,\ldots,s,\qquad \alpha=1,\ldots,t,
\end{gather}
with $\det(\partial f^\alpha/\partial y^\beta)\neq 0$.  Thus, the $x^i$ are invariants of the action, whose common level sets prescribe the $t$-dimensional pseudo-group orbits in $M$.    Let
\begin{gather}\label{cartan determining equations}
X=x,\qquad F^{(n_\star)}\big(x,y,Y^{(n_\star)}\big)=0,
\end{gather}
be the involutive determining equations.  We note that $n_\star\geq n^\star$, the order of the pseudo-group, as the minimal order determining equations might need to be completed to involution, \cite{S-2002}.

For any $1 \leq n \leq \infty $, let $\mathcal C\n$ denote the {\it contact system} on $\Dn$ spanned by the contact forms $\Upsilon_A^b$ of order $0 \leq \#A < n$,~\eqref{cf}. The f\/irst step in Cartan's procedure is to restrict the contact system $\mathcal C^{(n_\star)}$  to the subbundle $\G^{(n_\star)}\subset \D^{(n_\star)}$.  The aim is to recast the determining system \eqref{cartan determining equations} in terms of the Pfaf\/f\/ian system
\begin{gather}
 X^i-x^i=0, \qquad \qquad i=1,\ldots,s,\qquad\alpha=1,\ldots,t, \qquad 0\leq \# A \leq n_\star-1,\nonumber\\
 \Upsilon^{s+\alpha}_A|_{F^{(n_\star)}(x,y,Y^{(n_\star)})=0}=\left.\Pa{dY^\alpha_A-\sum_{b=1}^m Y^\alpha_{A,b}dz^b}\right|_{F^{(n_\star)}(x,y,Y^{(n_\star)})=0}=0.\label{cartan ideal}
\end{gather}
For $k\geq 1$, let $Y_{[k]}=(Y^1_{[k]},\ldots,Y^{t_k}_{[k]})$ be  local parameterizations for the f\/ibers of the bundles $\G^{(k)} \to \G^{(k-1)}$, where $t_k=\dim \G^{(k)}-\dim \G^{(k-1)}$ is the f\/iber dimension.  The system \eqref{cartan ideal} is then equivalent to
\begin{gather}
X^i-x^i =0, \qquad i=1,\ldots,s,\nonumber\\
dY^\alpha-\sum_{a=1}^m L^\alpha_a(z,Y,Y_{[1]})dz^a =0, \qquad \alpha=1,\ldots,t,\nonumber\\
dY^j_{[1]}-\sum_{a=1}^m L^j_{[1],a}(z,Y,Y_{[1]},Y_{[2]})dz^a =0, \qquad j=1,\ldots,t_1,\nonumber\\
 \qquad \vdots\nonumber\\
dY^k_{[n_\star-1]}-\sum_{a=1}^m L^k_{[n_\star-1],a}(z,Y,Y_{[1]},\ldots,Y_{[n_\star]})dz^a =0, \qquad  k=1,\ldots,t_{n_\star-1},\label{cartan contact forms}
\end{gather}
where the functions $L^\alpha_a, \ldots, L^i_{[n_\star-1],a}$ are prescribed by the determining system \eqref{cartan determining equations}.  With the dif\/ferential forms \eqref{cartan contact forms} in hand, Cartan proceeds, in an inductive manner, to derive a system of invariant one-forms that serve to characterize the pseudo-group.

Since the forms $dY^1, \ldots, dY^t$ are right-invariant and the action of $\G$ on $\G^{(n_\star)}$ preserves the contact system $\mathcal C^{(n_\star)}$, the forms
\begin{gather*}
\omega^{s+\alpha}_{[0]}=\sum_{a=1}^m L^\alpha_a(z,Y,Y_{[1]})dz^a,\qquad \alpha=1,\ldots,t,
\end{gather*}
must likewise be right-invariant. These together with the invariant forms
\begin{gather}\label{const MC forms}
\omega^{i}_{[0]}=dx^i,\qquad i=1,\ldots,s,
\end{gather}
constitute a basis of horizontal forms, and hence $dz^1,\ldots, dz^m$ can be expressed as linear combinations of $\omega^1_{[0]},\ldots, \omega^m_{[0]}$. Their exterior derivatives have the form
\[
d\omega^b_{[0]}=\sum_{a=1}^m  d L^b_a(z,Y,Y_{[1]})\wedge dz^a=\sum_{a=1}^m  \omega^a_{[0]}\wedge\pi^b_a,\qquad b=1,\ldots,m,
\]
where the one-forms $\pi^b_a$ are certain linear combinations of $\omega^1_{[0]},\ldots,\omega^m_{[0]}$,
$dY^1,\ldots,dY^t$, and $dY^1_{[1]},\ldots, dY^{t_1}_{[1]}$.  The invariance of $\omega^1_{[0]},\ldots,\omega^m_{[0]}$ implies that
\[
\sum_{a=1}^m  \omega^a_{[0]}\wedge(R^\star_\psi(\pi^b_a)-\pi^b_a)=0,\qquad b=1,\ldots,m,
\]
for all $\psi \in \G$ such that the pull-back $R^\star_\psi(\pi^b_a)$ is def\/ined.  This means that
\[
R^\star_{\psi}(\pi^b_a)\equiv \pi^b_a\quad \text{mod } \omega^1_{[0]},\ldots,\omega^m_{[0]}.
\]
By the assumptions, $t_1=\dim \G^{(1)}-\dim \G^{(0)}$ of the $\pi^b_a$ are linearly independent modulo $\omega^1_{[0]},\ldots$, $\omega^m_{[0]},dY^1,\ldots,dY^t$.  Those $t_1$ dif\/ferential forms can be written as
\[
\pi^i\equiv \sum_{j=1}^{t_1} c^i_j dY^j_{[1]}+\sum_{\alpha=1}^q e^i_\alpha dY^\alpha\quad\text{mod }\omega^1_{[0]},\ldots,\omega^m_{[0]},\qquad i=1,\ldots,t_1,
\]
with $\det (c^i_j)\neq 0$.  The coef\/f\/icients $c^i_j$ and $e^i_\alpha$ may depend on the variables $z$, $Y$, and $Y_{[1]}$.  By adding suitable multiples of the $\omega^a_{[0]}$ we can write
\[
\pi^i\equiv  \omega^i_{[1]} \mod \omega^1_{[0]},\ldots,\omega^m_{[0]}, \qquad i=1,\ldots,t_1,
\]
where
\begin{gather}\label{cartan first order invariant forms}
\omega^i_{[1]}:= \sum_{j=1}^{t_1} c^i_j \left(dY^j_{[1]}-\sum_{b=1}^m L^j_{[1],b}(z,Y,Y_{[1]}) dz^b\right)+\sum_{\alpha=1}^q e^i_\alpha \left(dY^\alpha-\omega^{p+\alpha}_{[0]}\right).
\end{gather}
 Cartan, \cite[pp.~597--600]{C-1904}, now concludes that the one-forms $\omega^1_{[1]} , \ldots \omega^{t_1}_{[1]}$, are right-invariant.  These \emph{first order Cartan forms} are equivalent to our f\/irst order Maurer--Cartan forms \eqref{mc forms} in the sense that
\[
\text{span}\big\{  \omega^1_{[1]} , \ldots, \omega^{t_1}_{[1]}\big\} = \text{span}\{  \mu^a_{Z^b}|_{L^{(n_\star)}(Z,\mu^{(n_\star)})=0} \}.
\]
Next by computing the exterior derivatives of the f\/irst order Cartan forms \eqref{cartan first order invariant forms} and repeating the above procedure, Cartan derives $t_2$ linearly independent invariant second order Cartan forms, and so on, up to order $n_\star-1$.

The $r_{n_\star-1}=m+t_1+t_2+\cdots+t_{n_\star-1}$ invariant one-forms so constructed are collectively denoted by $\omega^1,\omega^2,\ldots,\omega^{r_{n_\star-1}}$ without the subscripts. Their exterior derivatives can be written as
\begin{gather}\label{cartan structure eq}
d\omega^i=\sum_{1\leq j < k \leq r_{n_\star-1}} C^i_{jk}\, \omega^j\wedge\omega^k + \sum_{j=1}^{r_{n_\star-1}}\sum_{\beta=1}^{t_{n_\star}}  A^i_{j\beta}\, \omega^j\wedge\overline\pi^\beta,\qquad i=1,\ldots,r_{n_\star-1},
\end{gather}
where
\[
(\overline\pi^1,\ldots,\overline\pi^{ t_{n_\star}})\equiv \bpa{dY^1_{[n_\star]},\ldots,dY^{t_{n_\star}}_{[n_\star]}}\quad \text{mod }\omega^1,\ldots,\omega^{r_{n^{\star-1}}}.
\]
These constitute {\it Cartan's structure equations}.
If the pseudo-group acts intransitively, the structure coef\/f\/icients $C^i_{jk}$, $A^i_{j\beta}$ may depend on the invariants $x^1,\ldots,x^s$.

\section{Examples}

In this section we illustrate the two structure theories with a pair of elementary intransitive Lie pseudo-group actions.

\begin{Example}\label{intransitive pg ex}
Let $\mathcal{G}$ be the inf\/inite-dimensional Lie pseudo-group
\begin{gather*}
X=x,\qquad Y=y f(x)+\phi(x),\qquad Z=z(f(x))^x+\psi(x),
\end{gather*}
where $f, \phi, \psi \in C^{\infty}(\R)$ and $f(x) >0$.  This Lie pseudo-group was introduced by Cartan, \cite{C-1904-0}, as an example with an essential invariant.

We f\/irst construct the structure equations using Cartan's structure theory. The involutive determining system is
\begin{gather}
X=x,\qquad Y_z=0,\qquad Z_y=0,\qquad Z_z=(Y_y)^x,\qquad Y_{yy}=0,\nonumber\\
 Z_{zz}=0,\qquad Z_{zx}=(Y_y)^x\left(\log Y_y + \frac{x Y_{xy}}{Y_y}\right).\label{ex cartan def eq}
\end{gather}
Thus, the f\/ibers of the bundle $\pi^2_0:\G^{(2)}\to \G^{(0)}$ are parameterized by
\[
(Y_x,Y_y,Z_x,Y_{xx},Y_{xy},Z_{xx}),
\]
and the determining system \eqref{ex cartan def eq} is equivalent to the Pfaf\/f\/ian system
\begin{gather*}
X-x=0,\\
\Upsilon^y|_{\mathcal{G}^{(2)}}=dY-Y_x\, dx - Y_y\, dy=0,\qquad\quad\;
\Upsilon^z|_{\mathcal{G}^{(2)}}= dZ-Z_x\, dx - (Y_y)^x\, dz = 0,\\
\Upsilon^y_x|_{\mathcal{G}^{(2)}} = dY_x-Y_{xx}\, dx - Y_{xy} \, dy = 0,\qquad
\Upsilon^y_y|_{\mathcal{G}^{(2)}} = dY_y-Y_{yx}\, dx=0,\\
\Upsilon^z_x|_{\mathcal{G}^{(2)}} = dZ_x-Z_{xx}\, dx - (Y_y)^x\left(\log Y_y + \frac{x Y_{xy}}{Y_y}\right)\, dz=0.
\end{gather*}
%
Cartan's algorithm yields the six invariant one-forms
\begin{gather*}
\omega^1=dx,\qquad
\omega^2=Y_x\, dx + Y_y\, dy,
\qquad \omega^3=Z_x\, dx + (Y_y)^x\, dz,\\
\omega^4=\mu^y_X|_{\mathcal{G}^{(2)}}=\frac{1}{Y_y}\,\Upsilon^y_y, \qquad
\omega^5=\mu^y_Y|_{\mathcal{G}^{(2)}}=\Upsilon^y_x-\frac{Y_x}{Y_y}\,\Upsilon^y_y,\\
\omega^6=\mu^z_X|_{\mathcal{G}^{(2)}}=\Upsilon^z_x-\frac{x Z_x}{Y_y}\,\Upsilon^y_y.
\end{gather*}
By computing their exterior derivatives we obtain Cartan's structure equations
\begin{gather}
d\omega^1=0,\qquad
d\omega^2=\omega^4\wedge\omega^1 + \omega^5 \wedge \omega^2,\qquad
d\omega^3=\omega^6\wedge \omega^1+ x\, \omega^5 \wedge \omega^3,\nonumber\\
d\omega^4=\omega^1\wedge\overline{\pi}^1+\omega^2\wedge\overline{\pi}^2+\omega^5\wedge\omega^4,\qquad
d\omega^5=\omega^1\wedge\overline{\pi}^2,\nonumber\\
d\omega^6=\omega^1\wedge\overline{\pi}^3+\omega^3\wedge(\omega^5+x\,\overline{\pi}^2)+x\,\omega^5\wedge\omega^6,
\label{ex Cartan structure eq}
\end{gather}
where
\begin{gather*}
\overline{\pi}^1=\mu^y_{XX}|_{\mathcal{G}^{(2)}},\qquad \overline{\pi}^2=\mu^y_{XY}|_{\mathcal{G}^{(2)}},\qquad \overline{\pi}^3=\mu^z_{XX}|_{\mathcal{G}^{(2)}}.
\end{gather*}

On the other hand, the computation of the Maurer--Cartan structure equations by the algorithm presented in Section \ref{structure pseudo-group section} proceeds as follows.  The inf\/initesimal generators
\[
\vv=\xi(x,y,z)\,\partial_x + \eta(x,y,z)\,\partial_y + \zeta(x,y,z)\,\partial_z=\br{\alpha (x)\, y + \beta (x)}\,\partial_y + \br{\alpha (x)\,x\,z + \gamma(x)}\,\partial_z
\]
of the pseudo-group action \eqref{ex cartan def eq} are the solutions of the inf\/initesimal determining system
\begin{gather}\label{ex inf det eq}
\xi=0,\qquad \eta_z=0,\qquad \zeta_y=0,\qquad \zeta_z=x\,\eta_y,
\end{gather}
which can be obtained by linearizing \eqref{ex cartan def eq} at the identity.
As in \eqref{mc lin relations}, the lift of \eqref{ex inf det eq} produces the linear relations
\begin{gather}\label{ex1 lin rels}
\mu^x=0,\qquad \mu^y_Z=0,\qquad \mu^z_Y=0,\qquad \mu^z_Z=X\mu^y_Y,
\end{gather}
among the f\/irst order Maurer--Cartan forms.  On account of \eqref{ex1 lin rels} and its f\/irst prolongation, the structure equations for the zero-th and f\/irst order Maurer--Cartan forms are
\begin{gather}
d\mu^y=\mu^y_X \wedge\mu^x+\mu^y_Y\wedge\mu^y+\mu^y_Z\wedge\mu^z=\mu^y_Y\wedge\mu^y,\nonumber\\
d\mu^z=\mu^z_X \wedge\mu^x+\mu^z_Y\wedge\mu^y+\mu^z_Z\wedge\mu^z=X\mu^y_Y\wedge\mu^z,\nonumber\\
d\mu^y_X=\mu^y_Y\wedge\mu^y_X+\mu^y_{XY}\wedge\mu^y,\label{ex OP structure eq}\\
d\mu^y_Y=0,\nonumber\\
d\mu^z_X=X\mu^y_Y\wedge\mu^z_X+(\mu^y_Y+X\mu^y_{XY})\wedge\mu^z.\nonumber
\end{gather}
The two sets of structure equations \eqref{ex Cartan structure eq} and \eqref{ex OP structure eq} are isomorphic provided we set $x = X$ and
$
\omega^1=0
$
in Cartan's structure equations \eqref{ex Cartan structure eq}.
\end{Example}

\begin{Example}\label{ex intransitive group}
As the second example we consider the action of a one-dimensional Lie group on~$\mathbb{R}^2$ by translations
\begin{gather}\label{intrans g ex}
X=x\neq 0,\qquad Y=y+a\,x,\qquad a\in\mathbb{R}.
\end{gather}
%
Cartan computes the structure equations for this group, \cite[p. 1345]{C-1937}, and f\/inds
\begin{gather}\label{cartan struc eq int group}
d\omega^1=0,\qquad d\omega^2=\frac{1}{x}\,\omega^1\wedge\omega^2,
\end{gather}
where
\[
\omega^1=dx,\qquad  \omega^2=dy-\frac{y}{x}\,dx.
\]
Equations \eqref{cartan struc eq int group} involve two independent invariant one-forms and a non-vanishing, variable structure coef\/f\/icient.  They obviously do not conform with the structure equations for a one-dimen\-sio\-nal abelian Lie group.

On the other hand, the Maurer--Cartan structure equations \eqref{pseudo-group structure eq} for the pseudo-group \eqref{intrans g ex} have the desired form. First, the minimal order involutive determining system for the group action \eqref{intrans g ex} is
\begin{gather}\label{Ex2 DE}
X=x,\qquad Y-y=x Y_x,\qquad Y_y=1.
\end{gather}
Linearization of \eqref{Ex2 DE} yields the inf\/initesimal determining equations
\begin{gather}\label{Ex2 IDE}
\xi=0,\qquad \eta=x\,\eta_x,\qquad \eta_y=0,
\end{gather}
for the inf\/initesimal generators $\vv = \xi (x,y)\, \partial _x + \eta (x,y)\, \partial _y$.  The lift of \eqref{Ex2 IDE} produces the linear relations
\begin{gather}\label{lifted eq g ex}
\mu^x=0,\qquad \mu^y=X\mu^y_X,\qquad \mu^y_Y=0,
\end{gather}
among the f\/irst order Maurer--Cartan forms.  It follows from \eqref{lifted eq g ex} that $\mu^y$ forms a basis for the Maurer--Cartan forms.  Its exterior derivative is given by
\begin{gather}\label{op struc eq int group}
d\mu^y=\mu^y_Y\wedge\mu^y=0,
\end{gather}
which agrees with the structure equation for a one-dimensional abelian Lie group.  As in Examp\-le~\ref{intransitive pg ex}, Cartan's structure equations~\eqref{cartan struc eq int group} become equivalent with~\eqref{op struc eq int group} once we set $\omega^1=0$.
\end{Example}

Since there is no abstract object to represent a pseudo-group, saying when two pseudo-group actions come from the ``same pseudo-group'' is more tricky than in the f\/inite-dimensio\-nal case of Lie group actions. The following def\/inition encapsulates Cartan and Vessiot's notion of isomorphism, \cite{C-1904,V-1903}.

\begin{Definition}\label{isomorphic pseudo-groups}
Two pseudo-group actions $\G_1$, $\G_2$ on manifolds $M_1$, $M_2$ are \emph{isomorphic}, written $\G_1 \sim \G_2$, if they have a common isomorphic \emph{prolongation}, meaning a pseudo-group $\G$ acting on a~manifold $M$, and surjective submersions $\pi_i \colon M \to M_i$, $i=1,2$, such that, for each $i=1,2$, there is a one-to-one correspondence between elements $\varphi \in \G$ and $\varphi_i \in \G_i$ satisfying $\pi_i \comp \varphi = \varphi_i \comp \pi_i$.
\end{Definition}

For example, two actions of the same f\/inite-dimensional Lie group are isomorphic, as one can take $M = M_1 \times M_2$ with the Cartesian product action.  Proof of the transitive property of isomorphisms, i.e.~$\G_1 \sim \G_2$ and $\G_2 \sim \G_3$, then $\G_1 \sim \G_3$, can be found in  Stormark,~\cite{S-2000}.

On the other hand,  two isomorphic pseudo-groups need not have the same Cartan structure equations. A basic illustration of this fact is provided by Example \ref{ex intransitive group}.  Clearly, the group action~\eqref{intrans g ex} is isomorphic to the group of translations of $\mathbb{R}$,
\begin{gather}\label{transitive translation}
Y=y+a,\qquad a\in \mathbb{R},
\end{gather}
which is characterized by the single invariant one-form~$\omega=dy$.  The Cartan structure equation of the latter action is, of course, $d\omega=0$, which obviously is not isomorphic to the structure equations \eqref{cartan struc eq int group}.  On the other hand, the Maurer--Cartan structure equation of the group action~\eqref{transitive translation} is again given by~\eqref{op struc eq int group}.  In fact, it can be proved,~\cite{F-2008,F-2009}, that isomorphic pseudo-group actions always possess isomorphic Maurer--Cartan equations.

The two examples above show that, when dealing with intransitive Lie pseudo-group actions, the Maurer--Cartan structure equations \eqref{pseudo-group structure eq} and Cartan's structure equations \eqref{cartan structure eq} do not agree. We refer the reader to~\cite{F-2008} for more examples.  The discrepancy between the two sets of structure equations is due to the inclusion of the horizontal forms $\omega^1_{[0]},\ldots,\omega^s_{[0]}$, cf.~\eqref{const MC forms},
in Cartan's version.  They do not appear in the Maurer--Cartan structure equations \eqref{pseudo-group structure eq} since, for a Lie pseudo-group action of the form \eqref{cartan pg}, the f\/irst $s$ zero-th order Maurer--Cartan forms vanish:
\[
\mu^i=0,\qquad i=1,\ldots,s.
\]
Restricting to a target f\/iber yields
\[
\omega^i_{[0]}=\sigma^i=-\mu^i=0,\qquad i=1,\ldots,s.
\]
On the other hand, for transitive Lie pseudo-group actions, the two sets of structure equations are equivalent since the relations \eqref{key equality} provide a one-to-one correspondence between the zero-th order Maurer--Cartan forms $\mu^1,\ldots,\mu^m$ and the invariant horizontal forms $\sigma^1,\ldots,\sigma^m$.

\section{Duality}

In this f\/inal section, we investigate the relationship between pseudo-group structure equations and the commutator relations among their inf\/initesimal generators.  As we will see, the Maurer--Cartan structure equations of Theorem \ref{MCseqs} are naturally dual to the commutator relations among the inf\/initesimal generators, in the same sense as the f\/inite-dimensional version in \eqref{g}, \eqref{MC}.

Under the identif\/ication of inf\/inite jets of local vector f\/ields \eqref{infinitesimal generator} with their Taylor expansions
\[
j_\infty \vv|_{z_0}\simeq \sum_{a=1}^m\sum_{\# A \geq 0} \zeta^a_A(z_0)  \frac{(z-z_0)^A}{A!}\pp{}{z^a}\,,
\]
the f\/iber $J^\infty TM|_{z_0}$ inherits a Lie algebra structure.   The monomial vector f\/ields
\[
\vv^A_a=  \frac{(z-z_0)^A}{A!}\pp{}{z^a},\qquad \#A\geq 0,\qquad a=1,\ldots,m,
\]
provide a basis for the vector space $J^\infty TM|_{z_0}$ and satisfy the Lie bracket relations
\begin{gather}\label{diffeomorphism commutators}
\com{\vv_a^A}{\vv_b^B}= {A,B\sminus a\choose A}\,\vv_b^{A,B\sminus a}-{B,A\sminus b\choose B}\,\vv_a^{B,A\sminus b}.
\end{gather}
In the above equation
\[
{A,B\sminus a\choose A} = \mcases{\frac{(A,B\sminus a)!}{A!\, (B\sminus a)!},& a \in B,\\0,& a \not \in B,}
\]
where $B\sminus a$ denotes the multi-index obtained by deleting one occurrence of $a$ from $B$. By direct inspection, we conclude that, as in the f\/inite-dimensional theory, the commutation relations~\eqref{diffeomorphism commutators} are directly dual to the Maurer--Cartan structure equations \eqref{diffeomorphism MCx}.

The duality between the Maurer--Cartan structure equations and the Lie brackets of jets of inf\/initesimal dif\/feomorphism generators extends straightforwardly to general Lie pseudo-group actions.

\begin{Theorem}
The Maurer--Cartan structure equations \eqref{pseudo-group structure eq} of a Lie pseudo-group $\G$ at a~target fiber $(\bti)^{-1}(Z)$ are dual to the Lie algebra structure equations for the fiber $\Ji \g_{|Z}$ of the jet bundle of its infinitesimal generators.
\end{Theorem}

The proof relies on the observation that the Lie algebra structure equations for $\Ji\g$ are obtained by imposing the constraints prescribed by the inf\/initesimal determining equations~\eqref{linearized det system} on equations \eqref{diffeomorphism commutators}, while the Maurer--Cartan structure equations of a Lie pseudo-group $\G\subset \D$ are, in turn, obtained by imposing the constraints dictated by the lifted version \eqref{mc lin relations} of the inf\/initesimal determining equations on \eqref{diffeomorphism MCx}.
The details can be found in \cite{F-2009}.

Finally, we note that the horizontal forms $\omega_{[0]}^i=dx^i$, $i=1,\ldots,s$, in \eqref{const MC forms}
are, naturally, invariant under the group of translations $X^i=x^i+a^i$.  Thus Cartan's equations \eqref{cartan structure eq} more appropriately ref\/lect the inf\/initesimal structure of the extended set of transformations
\begin{gather}\label{transitive transformation}
X^i=x^i+b^i, \qquad Y^\alpha=f^\alpha(x,y),\qquad i=1,\ldots,s,\qquad \alpha=1,\ldots,t.
\end{gather}
acting transitively on $M$. However, there is no guarantee that the transformations \eqref{transitive transformation} represent a Lie pseudo-group.  Indeed, for the Lie group action \eqref{intrans g ex} of Example \ref{ex intransitive group}, the extension \eqref{transitive transformation} has the form
\begin{gather*}
X=x+b \ne 0,\qquad Y=y+ax,\qquad a,b\in\mathbb{R},
\end{gather*}
which does not def\/ine a transformation group.

\subsection*{Acknowledgments}

We would like to thank Olle Stormark and anonymous referees for helpful
remarks and references that served to improve the paper.  The first author
was supported in part by NSF Grant DMS 08--07317; the second author by NSF
Grant OCE 06--21134; the third author by NSF Grant DMS 05--05293 and a
University of Minnesota Graduate School Doctoral Dissertation Fellowship.

\pdfbookmark[1]{References}{ref}
\LastPageEnding


\begin{thebibliography}{99}

\footnotesize\itemsep=0pt

\bibitem{A-1989}
 Anderson  I.M.,
The variational bicomplex, Utah State Technical Report, 1989,
available at \mbox{\url{http://www.math.usu.edu/~fg_mp/}}.

\bibitem{BCGGG-1991}
Bryant  R.L., Chern  S.S., Gardner  R.B., Goldschmidt  H.L., Grif\/f\/iths  P.A.,
Exterior dif\/ferential systems, {\it Mathematical Sciences Research Institute Publications}, Vol.~18, Springer-Verlag, New York, 1991.

\bibitem{C-1904-0}
Cartan  \'E., Sur la structure des groupes inf\/inis, in  Oeuvres Compl\`etes, Part~II, Vol.~2, Gauthier-Villars, Paris, 1953,  567--569.

\bibitem{C-1904}
Cartan  \'E., Sur la structure des groupes inf\/inis de transformations, in  Oeuvres Compl\`etes, Part~II, Vol.~2, Gauthier-Villars, Paris, 1953,  571--714.

\bibitem{C-1937}
Cartan  \'E., La structure des groupes inf\/inis, in  Oeuvres Compl\`etes, Part~II, Vol.~2, Gauthier-Villars, Paris, 1953, 1335--1384.


\bibitem{Chrastina}
Chrastina  J.,
The formal theory of dif\/ferential equations, Masaryk University, Brno,  1998.


\bibitem{Ehresmann}
Ehresmann C.,
Introduction \`a la th\'eorie des structures inf\/init\'esimales et des pseudo-groupes de Lie, in   G\'eometrie Dif\/f\'erentielle, Colloques Internationaux du Centre National de la Recherche Scientif\/ique, Strasbourg, 1953, 97--110.

\bibitem{K-1989}
Kamran N.,
Contributions to the study of the equivalence problem of \'Elie Cartan and its applications to partial and ordinary dif\/ferential equations, {\it Acad. Roy. Belg. Cl. Sci. Mem. Collect. 8$_{\rm o}$ (2)} {\bf 45} (1989), no.~7, 122~pages.

\bibitem{K-1959}
Kuranishi M.,
On the local theory of continuous inf\/inite pseudo-groups.~I,
{\it Nagoya Math. J.} {\bf 15} (1959), 225--260.

\bibitem{K-1961}
Kuranishi M.,
On the local theory of continuous inf\/inite pseudo-groups.~II,
{\it Nagoya Math. J.} {\bf 19} (1961), 55--91.

\bibitem{L-1883}
Lie S., \"Uber unendlichen kontinuierliche Gruppen,
{\it Christ. Forh. Aar.} {\bf 8} (1883), 1--47
(see also Gesammelte Abhandlungen, Vol.~5, B.G.~Teubner, Leipzig, 1924, 314--360).

\bibitem{L-1891}
Lie S., Die Grundlagen f\"ur die Theorie der unendlichen kontinuierlichen Transformationsgruppen, {\it Leipzig. Ber.} {\bf 43} (1891), 316--393 (see also Gesammelte Abhandlungen, Vol.~6, B.G.~Teubner, Leipzig, 1927, 300--364).

\bibitem{LisleReidi}
Lisle I.G., Reid G.J.,
Cartan structure of inf\/inite Lie pseudogroups, in  Geometric Approaches to Dif\/ferential Equations  (Canberra, 1995), Editors P.J.~Vassiliou and I.G.~Lisle, {\it Austral. Math. Soc. Lect. Ser.}, Vol.~15, Cambridge University Press, Cambridge, 2000, 116--145.

\bibitem{Mackenzie}
Mackenzie K.,
Lie groupoids and Lie algebroids in dif\/ferential geometry, {\it London Mathematical Society Lecture Note Series}, Vol.~124, Cambridge University Press, Cambridge, 1987.

\bibitem{O-1995}
Olver P.J.,
Equivalence, invariants, and symmetry, Cambridge University Press, Cambridge, 1995.

\bibitem{OP-2005}
Olver P.J., Pohjanpelto J.,
Maurer--Cartan equations and structure of Lie pseudo-groups,
{\it Selecta Math. (N.S.)} {\bf 11} (2005), 99--126.


\bibitem{OP-2008}
Olver P.J., Pohjanpelto J.,
Moving frames for Lie pseudo-groups,
{\it Canad. J. Math.} {\bf 60} (2008), 1336--1386.

\bibitem{OP-2009}
Olver P.J., Pohjanpelto J.,
Dif\/ferential invariant algebras of Lie pseudo-groups, {\it Adv. Math.}, to appear.

\bibitem{P-1978}
Pommaret J.-F.,
Systems of partial dif\/ferential equations and Lie pseudogroups, {\it Mathematics and Its Applications}, Vol.~14, Gordon \& Breach Science Publishers, New York, 1978.

\bibitem{S-2002}
Seiler W.M.,
Involution -- the formal theory of dif\/ferential equations and its applications in computer algebra and numerical analysis, Habilitation Thesis, Dept. of Mathematics,  Universit{\"a}t Mannheim, 2002.

\bibitem{SS-1965}
Singer I., Sternberg S.,
The inf\/inite groups of Lie and Cartan. I.~The transitive groups,
{\it J. Analyse Math.} {\bf 15} (1965), 1--114.

\bibitem{S-2000}
Stormark O., Lie's structural approach to PDE systems, {\it Encyclopedia of Mathematics and Its Applications}, Vol.~80, Cambridge University Press, Cambridge, 2000.

\bibitem{F-2008}
Valiquette F.,
Structure equations of Lie pseudo-groups,
{\it J. Lie Theory} {\bf 18} (2008), 869--895.

\bibitem{F-2009}
Valiquette F.,
Applications of moving frames to Lie pseudo-groups, Ph.D. Thesis, University of Minnesota, 2009.

\bibitem{V-1903}
Vessiot E.,
Sur la th\'eorie des groupes continues,
{\it Ann. \'Ecole Norm. Sup.} {\bf 20} (1903), 411--451.

\end{thebibliography}
\end{document}